\documentclass[11pt]{amsart}
\usepackage[margin=1.1in]{geometry}
\usepackage{array}
\usepackage[dvipsnames]{xcolor}
\usepackage{amsmath,amssymb,amsthm,mathtools,enumitem,booktabs}
\usepackage{xr-hyper}
\usepackage{hyperref}
\hypersetup{colorlinks=true,linkcolor=blue!50!black,citecolor=blue!50!black,urlcolor=blue!50!black}

\newtheorem{theorem}{Theorem}[section]
\newtheorem{lemma}[theorem]{Lemma}

\newtheorem{corollary}[theorem]{Corollary}
\newtheorem{conjecture}[theorem]{Conjecture}
\theoremstyle{definition}

\theoremstyle{remark}
\newtheorem{remark}[theorem]{Remark}

\newcommand{\R}{\mathbb{R}}
\newcommand{\Z}{\mathbb{Z}}
\newcommand{\gs}{g^\star}
\newcommand{\Pdir}{P_{\mathrm{dir}}}
\newcommand{\Psym}{P_{\mathrm{sym}}}
\newcommand{\gsym}{g_{\mathrm{sym}}}
\newcommand{\capp}{\operatorname{cap}}
\newcommand{\Lg}{\log\!\big(2+M/n\big)}
\DeclareMathOperator{\Var}{Var}
\renewcommand{\Pr}{\mathbb{P}}
\newcommand{\E}{\mathbb{E}}

\title[Graphs versus loopless digraphs with a given degree sequence]{Graphs versus loopless digraphs with a given degree sequence,\\ and convergence of uniformly random graphs at every density}
\author{Behzad Mehrdad}
\date{November 1, 2026}

\begin{document}

\begin{abstract}
For a degree sequence $d$ on $n$ vertices with $M$ edges let $G(d)$ be the number of simple graphs and $D(d)$ the number of loopless digraphs with in-degree $=$ out-degree $=d$. We conjecture that $G(d)\ge\sqrt{D(d)}\,e^{-Cn\log(2+M/n)}$ for every graphical $d$, a single inequality across all densities, and prove it up to an extra factor $\log n$: $G(d)\ge\sqrt{D(d)}\,e^{-Cn\log n\log(2+M/n)}$ for every graphical $d$. The proof compares pairs of graphs with Eulerian digraphs by alternately colouring closed trails; the failure of injectivity is measured exactly, shared edges cancel against $2$-cycles, and parity is handled by merging trails and by a parity argument on Eulerian circuits of odd components. No asymptotic enumeration and no bound on the maximum-entropy multipliers is used. When the maximum-entropy edge probabilities lie in $[\delta,1-\delta]$, the same method gives the conjecture itself, $G(d)\ge\sqrt{D(d)}\,e^{-Cn\log n/\delta}$, by transferring three large-deviation estimates from independent arcs to uniform digraphs.

Combined with an exact identity $G(d)=e^{\gs(d)}\pi(d)$, where $\pi(d)$ is the probability that the maximum-entropy $\beta$-model has degrees exactly $d$, and a van der Waerden--type bound for $D(d)$, this shows that uniformly random graphs with degree sequence $d$ are close in the normalized cut metric to their $\beta$-model whenever $M\gg n\log n\log(M/n)$, at every density, with no tameness and no cap on the degrees. We also record the conjecture in the sparse regime and, with constant $\frac12+o(1)$, in the dense tame regime, from the asymptotic enumerations of McKay--Wormald and Barvinok--Hartigan.
\end{abstract}

\maketitle

\section{Introduction}

Counting graphs with a given degree sequence is understood asymptotically in two separate regimes: sparse (McKay--Wormald \cite{McK85,MW91}) and dense (Barvinok--Hartigan \cite{BH13}). Bipartite counts, by contrast, satisfy a uniform van der Waerden--type lower bound from Gurvits' theory of stable polynomials \cite{Gur08,Bar10}. The obstruction for graphs is that $\prod_{i<j}(1+x_ix_j)$ is not real stable. This paper proposes a single inequality comparing graphs with loopless digraphs, proves it up to a factor $\log n$ in the exponent for every degree sequence, and uses it to show that uniformly random graphs with given degrees converge to the maximum-entropy $\beta$-model at every density. A volumetric bound of this kind was first conjectured, in the form $G(d)\ge e^{\gs(d)-\eta n\log n}$, in \cite[Conj.~2.18]{Meh14}, where it was shown to imply that local tree statistics of uniformly random graphs agree with those of the $\beta$-model. The sparse limit theory that motivated this work, based on Stein's method with swaps, is developed in \cite{Meh-llogl}.

\subsection{Setup}
Let $d=(d_1,\dots,d_n)$ be nonnegative integers and $M=\tfrac12\sum_i d_i$.
\begin{itemize}[leftmargin=2em]
\item $G(d)$ is the number of simple graphs on $[n]$ with degree sequence $d$.
\item $D(d)$ is the number of loopless digraphs on $[n]$ in which vertex $i$ has out-degree and in-degree $d_i$; equivalently, $0$-$1$ matrices with zero diagonal and all row and column sums equal to $d$.
\end{itemize}
Isolated vertices change neither $G$ nor $D$, so we assume $d_i\ge1$ for all $i$; then $M\ge n/2$.

\begin{conjecture}[Main conjecture]\label{conj:main}
There is an absolute constant $C$ such that for every graphical $d$ with all $d_i\ge1$,
\[
G(d)\;\ge\;\sqrt{D(d)}\,\exp\!\big(-C\,n\Lg\big).
\]
\end{conjecture}

In the dense case $M\asymp n^2$ the loss is $\exp(-O(n\log n))$; for bounded average degree it is $\exp(-O(n))$. The hypothesis ``graphical'' is necessary: $d=(1,1,1)$ has $D(d)=2$ (two directed triangles) but $G(d)=0$. Since $G(d)^2\le D(d)\,e^{O(n\Lg)}$ (Corollary~\ref{cor:upper} and Theorem~\ref{thm:Dtwo}), the conjecture says that $G(d)$ is the square root of $D(d)$ at this scale.

\subsection{Main results}
Let $\gs(d)$ be the maximum entropy of a random graph with independent edges and expected degrees $d$, and $\pi(d)$ the probability that this maximum-entropy $\beta$-model has degrees exactly $d$.

\begin{theorem}[Theorem~\ref{thm:euler}, Lemma~\ref{lem:oddcomp}, Corollary~\ref{cor:euler}]\label{thm:A}
For every graphical $d$,
\[
G(d)\ \ge\ \sqrt{D(d)}\,\exp\big(-C\,n\log n\,\Lg\big).
\]
\end{theorem}

\begin{theorem}[Theorem~\ref{thm:euler2}]\label{thm:A2}
If the maximum-entropy edge probabilities of $d$ lie in $[\delta,1-\delta]$ and $n\ge n_0(\delta)$, then $G(d)\ge\sqrt{D(d)}\,e^{-Cn\log n/\delta}$ with $C$ absolute. This is Conjecture~\ref{conj:main} for such $d$, with the constant $O(1/\delta)$, proved without asymptotic enumeration.
\end{theorem}

\begin{theorem}[Corollary~\ref{cor:euler}]\label{thm:B}
If $d=d^{(n)}$ is strictly graphical and $M/(n\log n\,\Lg)\to\infty$, then a uniformly random graph $G_d$ with degree sequence $d$ satisfies $\|W^{G_d}-W^{w^\star}\|_\square=o(M/n^2)$ in probability, where $w^\star$ are the maximum-entropy probabilities. This holds at every density, with no tameness assumption and no cap on the degrees.
\end{theorem}

The other results of the paper are summarized in the following table.

\begin{center}
\renewcommand{\arraystretch}{1.08}
{\small\begin{tabular}{@{}>{\raggedright\arraybackslash}p{0.3\textwidth}>{\raggedright\arraybackslash}p{0.62\textwidth}@{}}
\toprule
Result & Statement\\
\midrule
Thm.~\ref{thm:identity} & $G(d)=e^{\gs}\pi(d)$ exactly\\
Thm.~\ref{thm:bip}, Thm.~\ref{thm:Dtwo} & $0$-$1$ matrices with forbidden positions; $e^{2\gs}(7e)^{-n}\prod(1+d_i)^{-1}\le D(d)\le e^{2\gs}$\\
Thm.~\ref{thm:reduction} & Conj.~\ref{conj:main} $\Leftrightarrow$ $\pi(d)\ge e^{-Cn\log(2+M/n)}$\\
Thm.~\ref{thm:sparse}, Thm.~\ref{thm:tame} & Conj.~\ref{conj:main} in the sparse ($d_{\max}=o(M^{1/4})$) and dense tame regimes, via \cite{McK85,MW91,BH13}\\
Thm.~\ref{thm:euler}, Thm.~\ref{thm:euler2}, Cor.~\ref{cor:euler} & Euler colouring: Theorems~\ref{thm:A} and \ref{thm:A2}\\
Thm.~\ref{thm:degconv}, Cor.~\ref{cor:mainconv}, Thm.~\ref{thm:degsparse} & transfer of convergence; Theorem~\ref{thm:B}; convergence in McKay's range\\
\bottomrule
\end{tabular}}
\end{center}

\subsection{Companion papers}
The sparse limit theory that motivated this work is developed in \cite{Meh-llogl}: by Stein's method of exchangeable pairs with degree-preserving swaps, uniformly random graphs with degrees $d$ converge to $\varphi\otimes\varphi$ for every $L^1$ degree profile with $M/n\to\infty$ and $M/n^2\to0$, with no cap on $d_{\max}$; the degree classes form a regularity partition with $O(\Lambda/\varepsilon^2+\varepsilon^{-1}\log\frac1\varepsilon)$ parts; and the $\beta$-graphon equals $\varphi\otimes\varphi$ in the sparse sense. The paper \cite{Meh-switch} collects other routes to bounds on $\pi(d)$: smoothing, which gives the conjecture for a nearby sequence; elementary switching bounds that avoid Gurvits' theorem; a volume bound $\pi(d)\ge e^{-Cn}\prod_i(1+\sigma_i)^{-1}$ in the sparse range, which removes the extra logarithm of Theorem~\ref{thm:A} there; the swap dynamics at positive density; and a continuum $f$-factor theorem. Open problems and conjectures arising from this work, including a volume form of Conjecture~\ref{conj:main} for all strictly graphical $d$, the sharp form $G\ge cM^{-1/4}\sqrt D$, and a programme for sparse graph limits beyond $L^p$, are collected in \cite{Meh-open}.

\subsection{Relation to Chatterjee--Diaconis--Sly}
The maximum-entropy random graph $G_W$ used throughout is the \emph{$\beta$-model} of Chatterjee, Diaconis and Sly \cite{CDS11}, with $\beta_i=-\lambda_i$. Their results line up with ours as follows.
\begin{center}
\renewcommand{\arraystretch}{1.25}
\begin{tabular}{@{}p{0.47\textwidth}p{0.47\textwidth}@{}}
\toprule
Chatterjee--Diaconis--Sly \cite{CDS11} & This paper\\
\midrule
$\beta$-model $p_{ij}=e^{\beta_i+\beta_j}/(1+e^{\beta_i+\beta_j})$; conditioned on its degrees it is uniform & $G_W$; the exact identity, Theorem~\ref{thm:identity}, and the tilt identity, \cite[Lemma~\ref{S-lem:tilt}]{Meh-switch}\\
MLE equations $d_i=\sum_jp_{ij}$ & maximum-entropy duality, Lemma~\ref{lem:expfamily}\\
Thm.~1.4: the closure of the set of expected degree sequences is $\operatorname{conv}(\text{graphical})$ & the degree polytope $\mathcal D_n$; strictly graphical $\Leftrightarrow$ no frozen pairs (\cite[Lemma~\ref{S-lem:nofrozen}]{Meh-switch})\\
Prop.~1.2: continuum Erd\H{o}s--Gallai inequality, strict in the interior & strict graphicality, in the dense scaling limit\\
Lemma~4.1: an Erd\H{o}s--Gallai margin of order $n^2$, with all $d_i\asymp n$, forces $\|\hat\beta\|_\infty=O(1)$ & \cite[Remark~\ref{S-rem:cdsquant}]{Meh-switch}; a relative margin alone does not bound the multipliers (\cite[Example~\ref{S-ex:hubcore}]{Meh-switch})\\
Lemma~6.2: in the tame case, $\Pr[\deg=d]\ge\frac12\exp(-c_\delta n^{3/2+\varepsilon})$ & a weak form of the bound on $\pi(d)$; Barvinok--Hartigan give $n^{-O(n)}$, and Theorem~\ref{thm:reduction} asks for $e^{-Cn\log(2+M/n)}$ uniformly; compare \cite[Proposition~\ref{S-prop:abscost}]{Meh-switch}\\
\bottomrule
\end{tabular}
\end{center}


\section{Notation and tools}

\paragraph{Entropy.} For $p\in[0,1]$ let $H(p)=-p\log p-(1-p)\log(1-p)$. $H$ is strictly concave, $H'(p)=\log\frac{1-p}{p}$, and $H'(p)\to+\infty$ as $p\to0^+$, $H'(p)\to-\infty$ as $p\to1^-$. For $u>0$, $z\in[0,1]$:
\begin{equation}\label{eq:legendre}
H(z)+z\log u\;\le\;\log(1+u),\qquad\text{with equality iff } z=\tfrac{u}{1+u}.
\end{equation}

\paragraph{Polytopes.} $\Pdir(d)$ is the set of real $n\times n$ matrices $Z$ with $0\le z_{ij}\le1$, $z_{ii}=0$, and all row and column sums equal to $d$; $\Psym(d)$ is its subset of symmetric matrices. Set
\[
g(Z)=\sum_{i\ne j}H(z_{ij}),\qquad \gsym(Z)=\sum_{i<j}H(z_{ij}),\qquad \gs(d)=\max_{Z\in\Psym(d)}\gsym(Z).
\]
When $d$ is graphical, $\Psym(d)\neq\emptyset$.

\paragraph{Forbidden positions.} For $r\in\Z_{\ge0}^m$, $c\in\Z_{\ge0}^n$ and $E\subseteq[m]\times[n]$, let $N(E;r,c)$ be the number of $0$-$1$ matrices supported on $E$ with row sums $r$ and column sums $c$; $P(E;r,c)$ the polytope of real matrices with entries in $[0,1]$, zero off $E$, and the same margins; and
\[
h(E;r,c)=\max_{Z\in P(E;r,c)}\sum_{(i,j)\in E}H(z_{ij}).
\]
Note $D(d)=N(E;d,d)$ with $E=\{(i,j):i\ne j\}$.

\paragraph{Real stability.} $p\in\R[x_1,\dots,x_k]$ is \emph{real stable} if $p(x)\ne0$ whenever $\operatorname{Im}x_i>0$ for all $i$. Products of real stable polynomials and elementary symmetric polynomials are real stable, and substituting $x_i=\ell_i(u)$ with $\ell_i$ nonzero linear forms with nonnegative coefficients preserves real stability.

\begin{theorem}[Gurvits \cite{Gur08}]\label{thm:gurvits}
Let $p$ be real stable, homogeneous of degree $K$ in $K$ variables, with nonnegative coefficients, and $\capp(p)=\inf_{u>0}p(u)/(u_1\cdots u_K)$. Then
\[
\frac{\partial^Kp}{\partial u_1\cdots\partial u_K}\;\ge\;\frac{K!}{K^K}\,\capp(p).
\]
\end{theorem}

\begin{theorem}[Darroch \cite{Dar64}]\label{thm:darroch}
Let $S$ be a sum of independent Bernoulli variables with $\E S=\mu\in\Z$. Then $\mu$ is the unique mode of $S$.
\end{theorem}

\paragraph{Stirling bounds.} For $k\ge1$: $(k/e)^k\sqrt{2\pi k}\le k!\le e\,k^{k+1/2}e^{-k}$.

\section{Symmetrization and the entropy upper bound}

\begin{lemma}[Symmetrization]\label{lem:sym}
If $\Pdir(d)\ne\emptyset$, then $g$ has a unique maximizer $Z^\star$ on $\Pdir(d)$, $Z^\star$ is symmetric, and $\max_{\Pdir(d)}g=2\gs(d)$. Moreover $\Psym(d)\neq\emptyset$ iff $\Pdir(d)\ne\emptyset$.
\end{lemma}
\begin{proof}
$\Pdir(d)$ is compact and convex and $g$ is strictly concave in the off-diagonal coordinates, so $Z^\star$ is unique. Transposition maps $\Pdir(d)$ to itself and preserves $g$, so $(Z^\star)^{\mathsf T}$ is also a maximizer, whence $Z^\star$ is symmetric. For symmetric $Z$, $g(Z)=2\gsym(Z)$, so $\max g=2\gsym(Z^\star)\le2\gs(d)$; conversely $2\gsym(Y)=g(Y)\le g(Z^\star)$ for every $Y\in\Psym(d)$. Finally $(Z+Z^{\mathsf T})/2\in\Psym(d)$ for $Z\in\Pdir(d)$.
\end{proof}

\begin{lemma}[Entropy bound]\label{lem:entropy}
For finite nonempty $S\subseteq\{0,1\}^N$,
$\displaystyle\log|S|\le\max_{z\in\operatorname{conv}(S)}\sum_{k=1}^NH(z_k).$
\end{lemma}
\begin{proof}
For a discrete random variable $Y$ write $\operatorname{Ent}(Y)=-\sum_y\Pr[Y=y]\log\Pr[Y=y]$ for its Shannon entropy. Let $X=(X_1,\dots,X_N)$ be uniform on $S$. Each point of $S$ has probability $1/|S|$, so
\[
\operatorname{Ent}(X)=-\sum_{s\in S}\frac1{|S|}\log\frac1{|S|}=\log|S|.
\]
By subadditivity of Shannon entropy, $\operatorname{Ent}(X)\le\sum_k\operatorname{Ent}(X_k)$. Each $X_k$ is $\{0,1\}$-valued with $\Pr[X_k=1]=\E X_k$, so $\operatorname{Ent}(X_k)=H(\E X_k)$, the binary entropy. Hence $\log|S|\le\sum_kH(\E X_k)$, and $\E X\in\operatorname{conv}(S)$.
\end{proof}

\begin{corollary}\label{cor:upper}
$D(d)\le e^{2\gs(d)}$, $G(d)\le e^{\gs(d)}$, and $N(E;r,c)\le e^{h(E;r,c)}$.
\end{corollary}
\begin{proof}
Apply Lemma~\ref{lem:entropy} to loopless digraphs (convex hull in $\Pdir(d)$) with Lemma~\ref{lem:sym}, and to graphs encoded by their entries above the diagonal (convex hull in $\Psym(d)$).
\end{proof}

\section{An exact identity for graphs}

Let $d$ be graphical and $W$ the unique maximizer of $\gsym$ on $\Psym(d)$. Call a pair $\{i,j\}$ \emph{free} if some $Y\in\Psym(d)$ has $y_{ij}\in(0,1)$, and \emph{frozen} otherwise; let $J$ be the set of free pairs. By convexity, every point of $\Psym(d)$ takes the same value in $\{0,1\}$ on each frozen pair.

\begin{lemma}\label{lem:expfamily}
$w_{ij}\in(0,1)$ for every free pair, and there are $\lambda_1,\dots,\lambda_n\in\R$ with
\[
\log\frac{1-w_{ij}}{w_{ij}}=\lambda_i+\lambda_j\qquad\text{for all }\{i,j\}\in J.
\]
\end{lemma}
\begin{proof}
For each free pair $p$ choose $Y^p\in\Psym(d)$ with $y^p_p\in(0,1)$ and let $Y^0$ be their average; every free coordinate of $Y^0$ lies in $(0,1)$. If a free coordinate of $W$ were $0$ or $1$, the one-sided derivative of $\gsym$ along $W+t(Y^0-W)$ at $t=0^+$ would be $+\infty$ (by the behaviour of $H'$ at the endpoints), contradicting maximality.

View $\gsym$ as a function of the free coordinates $z_J$. The feasible set is $\{z_J:Az_J=b\}\cap[0,1]^J$, with $A$ the vertex--pair incidence matrix restricted to $J$. Since $Y^0$ is interior to the box, the affine hull of the feasible set is $\{Az_J=b\}$. As $W$ is a relative-interior maximizer of a differentiable concave function, $\nabla\gsym(W)\perp\ker A$, so $\nabla\gsym(W)=\lambda^{\mathsf T}A$ for some $\lambda$. The $\{i,j\}$ entry of $\lambda^{\mathsf T}A$ is $\lambda_i+\lambda_j$, and $\partial\gsym/\partial z_{ij}=\log\frac{1-z_{ij}}{z_{ij}}$.
\end{proof}

Thus $w_{ij}=x_ix_j/(1+x_ix_j)$ on free pairs, with $x_i=e^{-\lambda_i}$. Let $G_W$ be the random graph in which each pair $\{i,j\}$ is an edge independently with probability $w_{ij}$, and
\[
\pi(d)=\Pr\big[\deg(G_W)=d\big].
\]

\begin{theorem}\label{thm:identity}
For every graphical $d$ and every graph $\Gamma$ with degree sequence $d$, $\Pr[G_W=\Gamma]=e^{-\gs(d)}$. Consequently
\[
G(d)=e^{\gs(d)}\,\pi(d).
\]
\end{theorem}
\begin{proof}
The adjacency matrix $\gamma$ of $\Gamma$ lies in $\Psym(d)$, so it agrees with $W$ on frozen pairs, which contribute factor $1$. On free pairs, by Lemma~\ref{lem:expfamily},
\[
\log\Pr[G_W=\Gamma]=\sum_J\log(1-w_{ij})-\sum_J\gamma_{ij}(\lambda_i+\lambda_j)
=\sum_J\log(1-w_{ij})-\sum_i\lambda_i(d_i-f_i),
\]
where $f_i$ is the number of frozen edges at $i$, the same for all points of $\Psym(d)$. The same computation with $W$ in place of $\gamma$ gives $\sum_Jw_{ij}(\lambda_i+\lambda_j)=\sum_i\lambda_i(d_i-f_i)$. Hence
\[
\log\Pr[G_W=\Gamma]=\sum_J\Big[\log(1-w_{ij})+w_{ij}\log\tfrac{w_{ij}}{1-w_{ij}}\Big]=-\sum_JH(w_{ij})=-\gs(d),
\]
since frozen pairs have $H=0$. Summing over $\Gamma$ gives $\pi(d)=G(d)e^{-\gs(d)}$.
\end{proof}

\begin{remark}\label{rem:dir}
By Lemma~\ref{lem:sym} the maximizer of $g$ on $\Pdir(d)$ is the same $W$; a pair is free for $\Pdir$ iff it is free for $\Psym$. The same argument with independent arcs gives $D(d)=e^{2\gs(d)}\pi_{\mathrm{dir}}(d)$, where $\pi_{\mathrm{dir}}(d)$ is the probability that independent arcs with probabilities $w_{ij}$ ($i\ne j$) produce all in- and out-degrees equal to $d$. Hence, exactly,
\[
\frac{G(d)}{\sqrt{D(d)}}=\frac{\pi(d)}{\sqrt{\pi_{\mathrm{dir}}(d)}}\;\ge\;\pi(d).
\]
\end{remark}

\section{A lower bound for 0-1 matrices with forbidden positions}

\begin{lemma}[One-variable extraction]\label{lem:onevar}
Let $p(t)=\prod_{k=1}^L(1+w_kt)$ with $w_k>0$ and $0\le r\le L$ an integer. Then
\[
[t^r]\,p\;\ge\;\frac{1}{7\sqrt{r+1}}\,\inf_{t>0}\frac{p(t)}{t^r},
\]
and the factor $1/(7\sqrt{r+1})$ may be replaced by $1$ when $r\in\{0,L\}$.
\end{lemma}
\begin{proof}
For $r=0$ both sides equal $1$; for $r=L$ both equal $\prod w_k$. For $0<r<L$, $p(t)/t^r\to\infty$ at both ends, so the infimum is attained at some $t_0>0$ with $t_0p'(t_0)/p(t_0)=r$. Let $S$ be a sum of independent Bernoulli$\big(w_kt_0/(1+w_kt_0)\big)$ variables. Then $\Pr[S=k]=[t^k]p\cdot t_0^k/p(t_0)$ and $\E S=r$, so $r$ is the mode of $S$ by Theorem~\ref{thm:darroch}. Let $\sigma^2=\Var S\le r$. If $\sigma=0$ then $\Pr[S=r]=1$. Otherwise Chebyshev gives $\Pr[|S-r|<2\sigma]\ge3/4$, and at most $4\sigma+1$ integers lie in that window, so
\[
\Pr[S=r]\ge\frac{3}{4(4\sigma+1)}\ge\frac{3}{16\sqrt r+4}\ge\frac{1}{7\sqrt{r+1}},
\]
the last since $21\sqrt{r+1}\ge16\sqrt r+4$. Rearranging gives the claim.
\end{proof}

\begin{lemma}[Weak duality]\label{lem:duality}
For all $x\in\R^m_{>0}$, $y\in\R^n_{>0}$ and $Z\in P(E;r,c)$,
\[
\frac{\prod_{(i,j)\in E}(1+x_iy_j)}{x^r\,y^c}\;\ge\;\exp\Big(\sum_{(i,j)\in E}H(z_{ij})\Big).
\]
In particular the infimum of the left side is at least $e^{h(E;r,c)}$.
\end{lemma}
\begin{proof}
Apply \eqref{eq:legendre} with $u=x_iy_j$ and sum: $\sum\log(1+x_iy_j)\ge\sum H(z_{ij})+\sum z_{ij}(\log x_i+\log y_j)=\sum H(z_{ij})+r\cdot\log x+c\cdot\log y$.
\end{proof}

\begin{theorem}\label{thm:bip}
Let $E\subseteq[m]\times[n]$, $r\in\Z^m_{\ge0}$, $c\in\Z^n_{\ge0}$ with $\sum r_i=\sum c_j=K$ and $P(E;r,c)\ne\emptyset$. Then
\[
N(E;r,c)\;\ge\;e^{h(E;r,c)}\prod_{i:\,r_i\ge1}\frac{1}{7\sqrt{r_i+1}}\prod_{j:\,c_j\ge1}\frac{1}{e\sqrt{c_j}}.
\]
\end{theorem}
\begin{proof}
Let $p(x,y)=\prod_{(i,j)\in E}(1+x_iy_j)$, so $N=[x^ry^c]p$, and let $E_i=\{j:(i,j)\in E\}$.

\emph{Step 1 (rows, exact).} $[x^r]p=q(y):=\prod_ie_{r_i}(y_{E_i})$, where $e_k$ is the $k$-th elementary symmetric polynomial. So $N=[y^c]q$, and $q$ is real stable, homogeneous of degree $K$, with nonnegative coefficients.

\emph{Step 2 (columns, Gurvits).} Discard columns with $c_j=0$ by setting $y_j=0$. Introduce variables $u_{j,1},\dots,u_{j,c_j}$ and set $Q(u)=q(\ell(u))$ with $\ell_j(u)=\frac1{c_j}\sum_ku_{j,k}$. Then $Q$ is real stable, homogeneous of degree $K$ in $K$ variables, with nonnegative coefficients. Only monomials of multidegree exactly $c$ survive the mixed derivative, and each $y_j^{c_j}$ contributes $c_j!/c_j^{c_j}$:
\[
\frac{\partial^KQ}{\prod_{j,k}\partial u_{j,k}}=[y^c]q\cdot\prod_j\frac{c_j!}{c_j^{c_j}}.
\]
By AM--GM, $\ell_j(u)\ge\tilde y_j:=(\prod_ku_{j,k})^{1/c_j}$, and $q$ is monotone in each variable, so $Q(u)/\prod u\ge q(\tilde y)/\tilde y^c$; taking $u_{j,k}=y_j$ gives equality. Hence $\capp(Q)=\inf_{y>0}q(y)/y^c$, and Theorem~\ref{thm:gurvits} with the Stirling bounds gives
\[
[y^c]q\;\ge\;\frac{K!}{\prod_jc_j!}\prod_j\Big(\frac{c_j}{K}\Big)^{c_j}\inf_{y>0}\frac{q(y)}{y^c}\;\ge\;\prod_{j:\,c_j\ge1}\frac{1}{e\sqrt{c_j}}\;\inf_{y>0}\frac{q(y)}{y^c}.
\]

\emph{Step 3 (back to rows, Darroch).} For fixed $y>0$, $e_{r_i}(y_{E_i})=[x_i^{r_i}]\prod_{j\in E_i}(1+x_iy_j)$, so by Lemma~\ref{lem:onevar}
\[
q(y)\;\ge\;\prod_{i:\,r_i\ge1}\frac{1}{7\sqrt{r_i+1}}\;\inf_{x>0}\frac{p(x,y)}{x^r}.
\]
Divide by $y^c$, take the infimum over $y$, and apply Lemma~\ref{lem:duality}.
\end{proof}

\begin{remark}
For $E=[m]\times[n]$ this is a quantitative form of Barvinok's estimate for $0$-$1$ contingency tables \cite{Bar10}. The loss is $\tfrac12\sum\log(1+r_i)+\tfrac12\sum\log(1+c_j)+O(m+n)$, which by Jensen is at most $\tfrac m2\log(1+K/m)+\tfrac n2\log(1+K/n)+O(m+n)$: it matches the $n\log(M/n)$ scale uniformly from sparse to dense. This is exactly what Conjecture~\ref{conj:main} asks for on the graph side.
\end{remark}

\section{Digraphs are determined by entropy; reduction of the conjecture}

\begin{theorem}\label{thm:Dtwo}
For every graphical $d$ with all $d_i\ge1$,
\[
e^{2\gs(d)}(7e)^{-n}\prod_{i=1}^n\frac{1}{1+d_i}\;\le\;D(d)\;\le\;e^{2\gs(d)}.
\]
In particular, with $c_0=7$,
$\;2\gs(d)-c_0\,n\Lg\le\log D(d)\le2\gs(d).$
\end{theorem}
\begin{proof}
The upper bound is Corollary~\ref{cor:upper}. For the lower bound apply Theorem~\ref{thm:bip} with $m=n$, $E=\{(i,j):i\ne j\}$, $r=c=d$: $P(E;d,d)=\Pdir(d)\ne\emptyset$ and $h(E;d,d)=2\gs(d)$ by Lemma~\ref{lem:sym}. The loss is at least $\prod_i(7\sqrt{d_i+1})^{-1}(e\sqrt{d_i})^{-1}\ge(7e)^{-n}\prod(1+d_i)^{-1}$. For the second form, Jensen gives $\sum\log(1+d_i)\le n\log(1+2M/n)\le n\log2+n\Lg$, and $\Lg\ge\log2$, so the loss is at most $(1+\log(14e)/\log2)\,n\Lg<7n\Lg$.
\end{proof}

\begin{theorem}[Reduction]\label{thm:reduction}
Let $C>0$.
\begin{enumerate}[leftmargin=2em]
\item If $\pi(d)\ge\exp(-Cn\Lg)$ for all graphical $d$ with $d_i\ge1$, then Conjecture~\ref{conj:main} holds with constant $C$.
\item If Conjecture~\ref{conj:main} holds with constant $C$, then $\pi(d)\ge\exp\!\big(-(C+c_0/2)\,n\Lg\big)$.
\end{enumerate}
\end{theorem}
\begin{proof}
(1) By Theorem~\ref{thm:identity} and Corollary~\ref{cor:upper}, $G(d)=e^{\gs}\pi(d)\ge\sqrt{D(d)}\,\pi(d)$.
(2) $\pi(d)=G(d)e^{-\gs}\ge\sqrt{D(d)}\,e^{-\gs}e^{-Cn\Lg}$, and $\sqrt{D(d)}\,e^{-\gs}\ge e^{-(c_0/2)n\Lg}$ by Theorem~\ref{thm:Dtwo}.
\end{proof}

Thus, up to the constant, the following are equivalent:
\begin{itemize}[leftmargin=2em]
\item $G(d)\ge\sqrt{D(d)}\,\exp(-Cn\Lg)$;
\item $G(d)\ge\exp(\gs(d)-Cn\Lg)$, the graph analogue of Theorem~\ref{thm:bip};
\item $\pi(d)\ge\exp(-Cn\Lg)$, anti-concentration of the $\beta$-model at its own mean.
\end{itemize}

\section{The sparse regime}\label{sec:sparse}

We use two classical formulas, valid as $M\to\infty$ with $d_{\max}=o(M^{1/4})$. Let
\[
\lambda=\frac{1}{4M}\sum_{i=1}^nd_i(d_i-1).
\]
\textbf{(MW)} Graphs \cite{McK85,MW91}:
\[
G(d)=\frac{(2M)!}{M!\,2^M\prod_id_i!}\exp\!\big(-\lambda-\lambda^2+o(1)\big).
\]
\textbf{(McK)} $0$-$1$ matrices with forbidden positions \cite{McK84}: for $S=\sum r_i=\sum c_j$ and a forbidden set $F$ with boundedly many entries per line,
\[
N=\frac{S!}{\prod r_i!\prod c_j!}\exp\!\Big(-\frac{\sum r_i(r_i-1)\sum c_j(c_j-1)}{2S^2}-\frac1S\sum_{(i,j)\in F}r_ic_j+o(1)\Big).
\]

\begin{theorem}\label{thm:sparse}
If $M\to\infty$ and $d_{\max}=o(M^{1/4})$, then
\[
G(d)=\big(1+o(1)\big)\,e^{1/2}(\pi M)^{-1/4}\sqrt{D(d)}.
\]
\end{theorem}
\begin{proof}
Apply (McK) with $r=c=d$, $S=2M$, $F$ the diagonal. The first exponent term is $(4M\lambda)^2/(2(2M)^2)=2\lambda^2$; the second is $\sum d_i^2/(2M)=2\lambda+1$. Hence
\[
\sqrt{D(d)}=\frac{\sqrt{(2M)!}}{\prod_id_i!}\exp\!\big(-\lambda^2-\lambda-\tfrac12+o(1)\big),
\qquad
\frac{G(d)}{\sqrt{D(d)}}=\frac{\sqrt{(2M)!}}{M!\,2^M}\,e^{1/2+o(1)}.
\]
By Stirling, $\sqrt{(2M)!}\sim(4\pi M)^{1/4}(2M/e)^M$ and $M!\,2^M\sim(2\pi M)^{1/2}(2M/e)^M$, with ratio $(\pi M)^{-1/4}$.
\end{proof}

\begin{corollary}
In the range of Theorem~\ref{thm:sparse}, Conjecture~\ref{conj:main} holds for every $C>0$ once $n$ is large, since $\log M\le2\log n=o(n)$.
\end{corollary}

\begin{remark}
The exact cancellation of $\lambda$ and $\lambda^2$ is the main evidence that $\sqrt{D(d)}$ is the right normalization: the degree-dependent corrections for graphs are precisely the square roots of those for loopless digraphs. The surviving $(\pi M)^{-1/4}$ comes from the single global parity constraint of graphs.
\end{remark}

\section{The dense tame regime}\label{sec:densetame}

Call $d$ \emph{$\delta$-tame} if $\delta\le w_{ij}\le1-\delta$ for all $i\ne j$; this forces $M\ge\delta n(n-1)/2$. Let $Q$ be the covariance matrix of $\deg(G_W)$:
\[
Q_{ii}=\sum_{j\ne i}w_{ij}(1-w_{ij}),\qquad Q_{ij}=w_{ij}(1-w_{ij})\ (i\ne j).
\]
\textbf{(BH)} Barvinok--Hartigan \cite{BH13}: for fixed $\delta>0$ and $\delta$-tame $d$, as $n\to\infty$,
\[
G(d)=\Theta_\delta(1)\,e^{\gs(d)}(2\pi)^{-n/2}(\det Q)^{-1/2}.
\]
(Their theorem gives the precise constant, including the parity factor and bounded corrections; only the order is used here.) By Theorem~\ref{thm:identity}, this is a local CLT: $\pi(d)=\Theta_\delta(1)(2\pi)^{-n/2}(\det Q)^{-1/2}$.

\begin{theorem}\label{thm:tame}
For fixed $\delta>0$ there are $c_\delta>0$ and $n_0(\delta)$ such that every $\delta$-tame $d$ with $n\ge n_0(\delta)$ satisfies
\[
\pi(d)\;\ge\;c_\delta\prod_{i=1}^n\frac{1}{\sqrt{2\pi d_i}}\;\ge\;c_\delta\exp\!\Big(-\frac n2\log\frac{4\pi M}{n}\Big).
\]
Consequently $G(d)\ge c_\delta\sqrt{D(d)}\exp(-\frac n2\log\frac{4\pi M}{n})$, i.e.\ Conjecture~\ref{conj:main} with constant $\tfrac12+o(1)$.
\end{theorem}
\begin{proof}
$Q$ is positive semidefinite, so by Hadamard's inequality $\det Q\le\prod_iQ_{ii}$, and $Q_{ii}\le\sum_jw_{ij}=d_i$. Insert into (BH); the second inequality is Jensen with $\sum d_i=2M$. The consequence follows from Theorem~\ref{thm:reduction}(1), since $M/n\to\infty$ makes $\log(4\pi M/n)/\log(2+M/n)\to1$.
\end{proof}

\begin{remark}
The same Hadamard step turns any local CLT for the $\beta$-model into the strong form, \cite[Conjecture~\ref{C-conj:strong}]{Meh-open}. The difficulty in the intermediate regime is therefore the local CLT itself, unavailable when $W$ has entries tending to $0$ or $1$ at different rates. Liebenau and Wormald \cite{LW24} obtained asymptotics across the gap for near-regular sequences; exact asymptotics for irregular sequences in the intermediate range remain open, while Corollary~\ref{cor:euler} gives Conjecture~\ref{conj:main} up to a factor $\log n$ in the exponent for every sequence.
\end{remark}

\section{Pairs of graphs and Eulerian digraphs}\label{sec:euler}

Alternately colouring the closed trails of a loopless digraph with in-degree $=$ out-degree $=d$ produces two graphs with degree sequence $d$. Every pair of graphs arises in this way, but not injectively. This section measures the failure of injectivity exactly and proves the main conjecture up to an extra factor $\log n$, for every graphical $d$.

Let $\mathcal D(d)$ be the set of loopless digraphs on $[n]$ with $\deg^+_i=\deg^-_i=d_i$ and no repeated arcs ($2$-cycles allowed), so $D(d)=|\mathcal D(d)|$. For $D\in\mathcal D(d)$ let $H(D)$ be the undirected graph of its $2$-cycles, $h_i=\deg_{H(D)}(i)$, and let $O(D)$ be the \emph{oriented part}: the arcs of $D$ not in a $2$-cycle. Then $O(D)$ has $\deg^\pm_i=d_i-h_i$ and no $2$-cycles. Call $D$ \emph{even} if every connected component of $O(D)$ has an even number of arcs, and let $D_{\rm even}(d)$ be the number of even $D$.

A \emph{transition system} on an Eulerian digraph $O$ is a choice, at every vertex $v$, of a bijection between the in-arcs and the out-arcs at $v$. It decomposes $O$ into closed directed trails; write $t(T)$ for their number. A transition system on a red--blue graph with equal red and blue degrees at every vertex is a choice, at every vertex, of a bijection between red and blue edges; it decomposes the graph into closed alternating trails, which have even length.

\begin{lemma}[Exact correspondence]\label{lem:eulerbij}
Fix an undirected graph $H$ with degrees $h\le d$. The following two sets are in bijection:
\begin{enumerate}[label=(\alph*)]
\item triples $(O,T,\chi)$: $O$ an oriented graph (no loops, repeated arcs or $2$-cycles) with $\deg^\pm=d-h$ and no arc on a pair of $H$; $T$ a transition system on $O$ all of whose trails have even length; $\chi$ an alternate colouring of the trails (two choices per trail);
\item triples $((R,B),T',o)$: $R,B$ edge-disjoint simple graphs avoiding $H$, each with degree sequence $d-h$; $T'$ an alternating transition system on $R\cup B$; $o$ an orientation of each of its trails.
\end{enumerate}
Moreover $D=O+\overleftrightarrow H$ ranges over the $D\in\mathcal D(d)$ with $H(D)=H$, and $(G_1,G_2)=(H\cup R,H\cup B)$ over the pairs in $\mathcal G(d)^2$ with $G_1\cap G_2=H$.
\end{lemma}
\begin{proof}
From (a): the red arcs of $O$ form a simple graph $R$, since two red arcs on one pair would be a $2$-cycle or a repeated arc of $O$; likewise $B$, and $R\cap B=\emptyset$ for the same reason. Consecutive arcs of a trail have different colours, so $T$ induces an alternating transition system $T'$, and $O$ induces orientations $o$. From (b): orienting the alternating trails gives an oriented graph, since $R\cap B=\emptyset$ and $R,B$ are simple; $T'$ becomes a directed transition system with even trails, and the colours are an alternate colouring. The two maps are inverse. For the last statement, $D=O+\overleftrightarrow H$ has $\deg^\pm=d$, its $2$-cycles are exactly $H$ because $O$ has none and avoids $H$, and $G_1\cap G_2=H$ because $R,B$ avoid $H$ and are disjoint.
\end{proof}

\begin{lemma}[Number of trails]\label{lem:trails}
Let $T$ be a uniformly random transition system, either on an Eulerian digraph with $\deg^\pm=c$ or on a red--blue graph with red and blue degrees $c$. Then
\[
E\,t(T)\le\sum_v\sum_{k=1}^{c_v}\frac1k\le\sum_v(1+\log(1+c_v)),\qquad E\,2^{t(T)}\le\prod_v(1+c_v).
\]
\end{lemma}
\begin{proof}
Reveal $T$ along an exploration of the trails: whenever the current trail enters $v$ and $v$ has $k$ unmatched outgoing arcs (blue edges), the partner is uniform among them. The trail closes only if the partner is the arc with which it started, which has conditional probability at most $1/k$. At $v$ the values of $k$ are $c_v,c_v-1,\dots,1$. So $t$ is a sum of indicators $X_s$ with $\Pr[X_s=1\mid\text{past}]\le1/k_s$, which gives both bounds, the second by induction on $s$ ($E[2^{X_s}\mid\text{past}]\le1+1/k_s$).
\end{proof}

\begin{lemma}[Parity]\label{lem:parity}
Let $O$ be an Eulerian digraph with $\deg^\pm=c$, $m=\sum_vc_v$ arcs and maximum degree $c_{\max}$, every component of which has an even number of arcs. For a uniformly random transition system $T$ on $O$,
\[
\Pr[\text{all trails of $T$ are even}]\ \ge\ \tfrac14\exp\Big(-2\log(2mc_{\max})\sum_v\sum_{k=1}^{c_v}\frac1k\Big).
\]
\end{lemma}
\begin{proof}
Within a component, any two trails are joined by a chain of trails, consecutive ones sharing a vertex; the number of odd trails in a component is even. Define $f(T)$ as follows: while odd trails exist, take (canonically) two odd trails $P,Q$ in one component at minimal distance in this chain graph, so that the trails $E_1,\dots,E_{j}$ in between are even, and merge $P$ with $E_1$, the result with $E_2$, \dots, the result with $Q$. Each merge swaps the transitions of two different trails at a common vertex $v$, $(e_1\to f_1),(e_2\to f_2)\mapsto(e_1\to f_2),(e_2\to f_1)$, which concatenates them; the final trail is even. So $f(T)$ has only even trails and is reached from $T$ by $r\le t(T)$ swaps. Each swap is undone by naming a vertex and two transitions at it, at most $\sum_v\binom{c_v}2\le mc_{\max}$ choices, so every $T'$ has at most $\sum_{r\le\tau}(mc_{\max})^r\le(2mc_{\max})^\tau$ preimages with $t(T)\le\tau$. With $\tau=2E\,t$ and Markov's inequality, $\Pr[t\le\tau]\ge\frac12$, and Lemma~\ref{lem:trails} gives the bound.
\end{proof}

\begin{theorem}[Euler colouring bound]\label{thm:euler}
For every graphical $d$ with $M\ge1$,
\[
G(d)^2\ \ge\ \frac{D_{\rm even}(d)}{\prod_v(1+d_v)}\cdot\frac14\exp\Big(-2\log(4Md_{\max})\sum_v\sum_{k=1}^{d_v}\frac1k\Big).
\]
In particular $G(d)\ge\sqrt{D_{\rm even}(d)}\,\exp\big(-Cn\log n\,\Lg\big)$.
\end{theorem}
\begin{proof}
Fix $H$ and count both sides of Lemma~\ref{lem:eulerbij}. Side (b) is $\sum_{(R,B)}\sum_{T'}2^{t(T')}\le\#\{(R,B)\}\prod_v(d_v-h_v)!\prod_v(1+d_v-h_v)$ by Lemma~\ref{lem:trails}. Side (a) is $\sum_O\sum_{T\text{ even}}2^{t(T)}\ge\sum_O\prod_v(d_v-h_v)!\,\Pr_T[\text{all even}]$. If $O+\overleftrightarrow H$ is even, Lemma~\ref{lem:parity} applies to $O$ with $m\le2M$ and $c\le d$. Divide by $\prod_v(d_v-h_v)!$, sum over $H$, and use $\prod(1+d_v-h_v)\le\prod(1+d_v)$. The last statement uses $\sum_vH_{d_v}\le n(1+\log(1+2M/n))$, $\log(4Md_{\max})\le3\log(2n)$ and $\prod_v(1+d_v)\le(1+2M/n)^n$.
\end{proof}

\begin{lemma}[Odd components]\label{lem:oddcomp}
For every $d$, $D_{\rm even}(d)\ge D(d)\,(8M^2)^{-n/6}$.
\end{lemma}
\begin{proof}
Let $C_1,C_2$ be two odd components of $O(D)$, with vertex sets $V_1,V_2$ and $L_1,L_2$ arcs. Since $O$ has no arcs between different components, every pair $\{x,y\}$ with $x\in V_1$, $y\in V_2$ is either \emph{full} (a $2$-cycle of $D$) or \emph{empty} (no arc). Call arcs $a\to b$ of $C_1$ and $c\to e$ of $C_2$ \emph{compatible} if the pairs $\{a,e\}$ and $\{b,c\}$ have the same type.

\emph{Compatible arcs exist.} Let $a_0\to a_1\to\dots\to a_{L_1}=a_0$ and $c_0\to\dots\to c_{L_2}=c_0$ be Eulerian circuits of $C_1$ and $C_2$, and let $x_{i,j}\in\{0,1\}$ be the type of $\{a_i,c_j\}$, with indices modulo $L_1$ and $L_2$. If no two arcs were compatible, then $x_{i,j+1}+x_{i+1,j}=1$ for all $i,j$ (take the arcs $a_i\to a_{i+1}$ and $c_j\to c_{j+1}$). Hence $x_{i+1,j}=1-x_{i,j+1}$, so moving from $(i,j)$ to $(i+1,j-1)$ flips the value. After $L_1L_2$ such moves one is back at $(i,j)$, and $L_1L_2$ is odd, a contradiction.

\emph{The switch.} Take compatible arcs $a\to b$, $c\to e$. If $\{a,e\},\{b,c\}$ are empty, replace $a\to b$, $c\to e$ by $a\to e$, $c\to b$. If they are full, replace $b\to c$, $e\to a$ by $b\to a$, $e\to c$; this is the same switch performed in the complement digraph, whose oriented part is the reverse of $O(D)$. In both cases all in- and out-degrees are preserved, no repeated arc is created, and the two new arcs replace $a\to b$, $c\to e$ in the oriented part: in the second case $\{a,b\},\{c,e\}$ become $2$-cycles and $\{a,e\},\{b,c\}$ become oriented. Other components are untouched. The new component is the closed trail $b\leadsto a\to e\leadsto c\to b$ through all arcs of $C_1\cup C_2$, so it is connected and has $L_1+L_2$ arcs, an even number.

\emph{Counting.} The number of arcs of $O(D)$ is $2M-2|H|$, so odd components come in pairs, and each has at least $3$ arcs and $3$ vertices, so there are at most $n/3$ of them. Repeating the switch canonically maps $\mathcal D(d)$ to even digraphs with at most $n/6$ switches. Each switch is undone by naming the two new oriented arcs and the type, at most $2(2M)^2$ choices.
\end{proof}

\begin{lemma}[Parity, sharpened]\label{lem:parity2}
Let $O$ be as in Lemma~\ref{lem:parity}, and let $V=\mathcal V\sqcup L$ be a partition of its vertex set, $K,\ell\ge1$, such that
\begin{enumerate}[label=(\roman*)]
\item any two vertices of $\mathcal V$ are joined by a path of length at most $K$ in $O$;
\item every vertex of $L$ in a component of $O$ meeting $\mathcal V$ is within distance $\ell$ of $\mathcal V$, and every component of $O$ not meeting $\mathcal V$ has diameter less than $\ell$.
\end{enumerate}
Then
\[
\Pr[\text{all trails of $T$ are even}]\ \ge\ \prod_v\frac1{\max(c_v,1)}\cdot(2mc_{\max})^{-(K+1)n/2-(2\ell+K+1)|L|}.
\]
With $L=\emptyset$ the hypothesis is that $O$ is connected with diameter at most $K$, and the bound is $\prod_v\max(c_v,1)^{-1}(2mc_{\max})^{-(K+1)n/2}$.
\end{lemma}
\begin{proof}
Reveal $T$ by the exploration of Lemma~\ref{lem:trails}. A trail can close only when it returns to its starting vertex $u$, with $k\ge1$ unmatched out-arcs there, one of which is its starting arc $s_0$. Call such a return \emph{bad} if closing now would give an odd trail. Let $A$ be the event that at every bad return with $k\ge2$ the partner is not $s_0$. Its conditional probability at each such step is $1-1/k$, and at a vertex $v$ the values of $k$ are distinct and run through $c_v,\dots,1$, so by sequential conditioning
\[
\Pr[A]\ge\prod_v\prod_{k=2}^{c_v}\Big(1-\frac1k\Big)=\prod_v\frac1{\max(c_v,1)}.
\]
On $A$ an odd trail can only be closed at a step with $k=1$, at its starting vertex, and a step with $k=1$ at a vertex is the last entry into that vertex, so it happens at most once there. Hence for $T\in A$ the odd trails have distinct starting vertices: at most $|\mathcal V|\le n$ of them start in $\mathcal V$ and at most $|L|$ start in $L$.

Apply the merging map of Lemma~\ref{lem:parity}, choosing each time two odd trails of one component at minimal distance in the chain graph of trails, so that the trails in between are even. A path of length $r$ in $O$ from a vertex of one trail to a vertex of another gives a chain of length at most $r+1$. Hence two odd trails started in $\mathcal V$ are at distance at most $K+1$, by (i); an odd trail started in $L$ is at distance at most $\ell+K+1$ from one started in $\mathcal V$ and at most $2\ell+K+1$ from another started in $L$ of the same component, by (ii) and (i), or at most $\ell$ if their component does not meet $\mathcal V$. The odd trails present at any stage are among the original ones, with the same starting vertices. A step that joins two trails started in $\mathcal V$ costs at most $K+1$ merges, and a step that consumes a trail started in $L$ costs at most $2\ell+K+1$ merges, since the pair chosen is at minimal distance. There are at most $n/2$ steps of the first kind and at most $|L|$ of the second, so at most $(K+1)n/2+(2\ell+K+1)|L|$ merges in all, and as in Lemma~\ref{lem:parity} each $T'$ with only even trails has at most $(2mc_{\max})^{(K+1)n/2+(2\ell+K+1)|L|}$ preimages in $A$.
\end{proof}

\begin{theorem}[Euler colouring bound at the dense scale]\label{thm:euler2}
There is an absolute constant $C$ such that for every $\delta\in(0,\frac12]$ there is $n_0(\delta)$ with the following property: every $\delta$-tame $d$ (Section~\ref{sec:densetame}) on $n\ge n_0(\delta)$ vertices satisfies
\[
G(d)\ \ge\ \sqrt{D(d)}\,\exp\Big(-\frac{C\,n\log n}{\delta}\Big).
\]
\end{theorem}
\begin{proof}
Let $w=w^\star$ be the maximum-entropy probabilities, $\delta\le w_{ij}\le1-\delta$. Under the model with independent arcs $i\to j$ of probabilities $w_{ij}$ every $D\in\mathcal D(d)$ has the same probability, $\pi_{\rm dir}(d)/D(d)=e^{-2\gs(d)}$ (Remark~\ref{rem:dir}), so the uniform measure on $\mathcal D(d)$ is this model conditioned on $\deg^\pm=d$, and for every property $\mathcal B$ of digraphs
\begin{equation}\label{eq:dirtransfer}
\frac{\#\{D\in\mathcal D(d):\mathcal B\}}{D(d)}\ \le\ \frac{\Pr_{\rm ind}[\mathcal B]}{\pi_{\rm dir}(d)}\ \le\ (7e)^n\prod_i(1+d_i)\Pr_{\rm ind}[\mathcal B]\ \le\ e^{3n\log n}\Pr_{\rm ind}[\mathcal B],
\end{equation}
by Theorem~\ref{thm:Dtwo}. Under the independent model, a pair $\{i,j\}$ carries exactly one arc, i.e.\ belongs to $O(D)$, with probability $q_{ij}=2w_{ij}(1-w_{ij})\in[p_0,\frac12]$, $p_0=2\delta(1-\delta)\ge\delta$, independently over pairs. Let $\varepsilon=p_0/4$, $\ell=\lceil160\,p_0^{-1}\log n\rceil$, let $\deg_O(v)$ be the number of pairs of $O=O(D)$ at $v$, $L=\{v:\deg_O(v)\le\varepsilon n\}$ and $\mathcal V=V\setminus L$. Consider three bad properties:
\begin{enumerate}[label=(B\arabic*)]
\item some disjoint $S,T$ with $|S|,|T|\ge\varepsilon n$ have fewer than $\frac{p_0}2|S||T|$ pairs of $O$ between them;
\item $|L|\ge n/\log n$;
\item $O$ contains a simple path with $\ell$ arcs all of whose vertices lie in $L$.
\end{enumerate}
We use Chernoff's bound $\Pr[X\le\mu/2]\le e^{-\mu/8}$ for a sum $X$ of independent indicators of mean $\mu$. For fixed disjoint $S,T$ the number of pairs of $O$ between them is such a sum with $\mu\ge p_0|S||T|$, so $\Pr_{\rm ind}[\text{(B1)}]\le4^ne^{-p_0\varepsilon^2n^2/8}$. For a set $U$ of $k\le n/2$ vertices, $U\subseteq L$ forces at most $k\varepsilon n=p_0kn/4$ pairs of $O$ between $U$ and $V\setminus U$, a sum of $k(n-k)\ge kn/2$ independent indicators with $\mu\ge p_0kn/2$; so $\Pr_{\rm ind}[U\subseteq L]\le e^{-p_0kn/16}$, and with $k=\lceil n/\log n\rceil$, $\Pr_{\rm ind}[\text{(B2)}]\le2^ne^{-p_0n^2/(16\log n)}$. For a fixed sequence of $\ell+1\le n/2$ distinct vertices with vertex set $U$, the $\ell$ path pairs lie inside $U$ and are independent of the pairs between $U$ and $V\setminus U$, so the sequence is a path inside $L$ with probability at most $2^{-\ell}e^{-p_0(\ell+1)n/16}$, and $\Pr_{\rm ind}[\text{(B3)}]\le n^{\ell+1}2^{-\ell}e^{-p_0(\ell+1)n/16}\le e^{-(\ell+1)p_0n/32}\le e^{-5n\log n}$ once $p_0n\ge32\log n$. For $n\ge n_0(\delta)$ all three are at most $\frac16e^{-4n\log n}$, so by \eqref{eq:dirtransfer} the bad digraphs number at most $\frac12e^{-n\log n}D(d)\le\frac12(8M^2)^{-n/6}D(d)$, and by Lemma~\ref{lem:oddcomp} at least $\frac12(8M^2)^{-n/6}D(d)$ digraphs in $\mathcal D(d)$ are even and good.

Let $D$ be even and good, $O=O(D)$. For $u,u'\in\mathcal V$, either $u'\in N(u)$, or $N(u)$ and $N(u')$ meet, or they are disjoint sets of size at least $\varepsilon n$, between which (B1) failing gives a pair of $O$, hence a path $u\,s\,t\,u'$; so (i) of Lemma~\ref{lem:parity2} holds with $K=3$. A shortest path from a vertex of $L$ to $\mathcal V$, or between two vertices of a component inside $L$, is a simple path all of whose vertices except possibly the last lie in $L$, so (B3) failing gives (ii) with this $\ell$. Lemma~\ref{lem:parity2}, with $|L|<n/\log n$, $m\le2M$, $c_{\max}\le d_{\max}$ and $2mc_{\max}\le2n^3$, gives
\[
\Pr_T[\text{all even}]\ \ge\ n^{-n}(2n^3)^{-2n-(2\ell+4)n/\log n}\ \ge\ \exp\Big(-\frac{C_1n\log n}{p_0}\Big).
\]
As in the proof of Theorem~\ref{thm:euler} (count both sides of Lemma~\ref{lem:eulerbij} for each $H$ and sum), $G(d)^2\ge\sum_{D\in\mathcal D(d)}\Pr_T[\text{all even}]\big/\prod_v(1+d_v)$, and restricting the sum to even good $D$ gives $G(d)^2\ge\tfrac12(8M^2)^{-n/6}D(d)\,e^{-C_1n\log n/p_0}n^{-n}$. Since $8M^2\le2n^4$ and $p_0\ge\delta$, the claim follows.
\end{proof}

\begin{remark}
For $\delta$-tame $d$ one has $M\ge\delta\binom n2$, so $n\Lg\asymp n\log n$, and Theorem~\ref{thm:euler2} is Conjecture~\ref{conj:main} for such $d$ with the constant $O(1/\delta)$, proved without asymptotic enumeration: besides switchings, its only external input is Gurvits' capacity bound through Theorem~\ref{thm:Dtwo}, used to transfer three large-deviation estimates for independent arcs to uniform digraphs. The transfer costs $e^{O(n\log n)}$, which is why the bad properties were chosen to have independent-model probability $e^{-\Omega(n^2/\log n)}$ or $e^{-5n\log n}$; the vertices of small degree in $O(D)$, which cannot be excluded at that scale, are handled by the partition of Lemma~\ref{lem:parity2}. Theorem~\ref{thm:tame} gives the constant $\frac12+o(1)$ via \cite{BH13}.
\end{remark}

\begin{corollary}\label{cor:euler}
For every graphical $d$,
\[
G(d)\ \ge\ \sqrt{D(d)}\,\exp\big(-Cn\log n\,\Lg\big),
\]
that is, Conjecture~\ref{conj:main} with an extra factor $\log n$. Consequently, for strictly graphical $d$, $-\log\pi(d)\le Cn\log n\,\Lg$, and by Theorem~\ref{thm:degconv} uniformly random graphs with degree sequence $d$ are cut-close to their $\beta$-model after normalization whenever $M/(n\log n\,\Lg)\to\infty$, at every density, with no tameness and no cap on the degrees.
\end{corollary}
\begin{proof}
Theorem~\ref{thm:euler} and Lemma~\ref{lem:oddcomp}. For $\pi$, use $G=e^{\gs}\pi$ and $D\ge e^{2\gs}(7e)^{-n}\prod(1+d_i)^{-1}$ (Theorem~\ref{thm:Dtwo}).
\end{proof}

\begin{remark}[Why this works where switchings do not]\label{rem:euler}
\begin{enumerate}[label=(\alph*)]
\item \emph{Shared edges cost nothing.} Edges in $G_1\cap G_2$ correspond exactly to $2$-cycles of $D$ and are put into $H$ on both sides of Lemma~\ref{lem:eulerbij}; the factor $2^{-\Theta(n^2)}$ that appears in the dense case if one insists on colouring $2$-cycles is never paid.
\item \emph{Multipliers never enter.} The argument compares counts directly and uses no $\beta$-model parameters, so it is insensitive to frozen pairs and to unbounded multipliers (\cite[Example~\ref{S-ex:hubcore}]{Meh-switch}), unlike the switching bounds of \cite[\S\ref{S-sec:switch}]{Meh-switch}.
\item \emph{The extra $\log n$} comes from paying $\log(4Md_{\max})$ per merge in Lemma~\ref{lem:parity}; Theorem~\ref{thm:euler2} removes it in the dense tame regime. Removing it in general is discussed in \cite{Meh-open}.
\end{enumerate}
\end{remark}

\section{Consequences: convergence of uniformly random graphs at every density}\label{sec:conv}

The lower bound on $\pi(d)$ is not only a counting statement: it controls the cost of passing from the product measure of a $\beta$-model to the uniform measure on graphs with degree sequence $d$ (\S\ref{subsec:transfer}), an idea going back to \cite{CDS11} in the dense case and to \cite{Meh14}. Combined with Corollary~\ref{cor:euler} it gives convergence to the $\beta$-model at every density with $M\gg n\log n\Lg$ (Corollary~\ref{cor:mainconv}). After the transfer inequality we state the corresponding theorem in McKay's range. Convergence without any bound on $\pi(d)$, by Stein's method with swaps, is proved for sparse degree sequences in \cite{Meh-llogl}, and \cite{Meh-switch} gives bounds on $\pi(d)$ by switchings.

Throughout, $G_d$ is uniform on graphs with degree sequence $d$, $\rho=2M/n^2$, $\bar d=2M/n$, $\varphi_n(x)=d_{\lceil nx\rceil}/\bar d$, and for a simple graph $G$ or a symmetric array $w$ with zero diagonal, $W^G$ and $W^w$ are the step graphons on the intervals $I_i=[\frac{i-1}n,\frac in)$. The cut norm is $\|U\|_\square=\sup_{S,T}|\int_{S\times T}U|$.

\subsection{Transfer}\label{subsec:transfer}
For $x\in(0,\infty)^n$ let $G_x$ be the $\beta$-model with edge probabilities $w_{ij}=x_ix_j/(1+x_ix_j)$, and $\pi_x(d)=\Pr[\deg(G_x)=d]$. As in Theorem~\ref{thm:identity}, every graph with degree sequence $d$ has the same probability under $G_x$, so for every event $\mathcal A$,
\begin{equation}\label{eq:transfer}
\Pr[G_d\in\mathcal A]\;\le\;\Pr[G_x\in\mathcal A]/\pi_x(d).
\end{equation}
For strictly graphical $d$ the maximum-entropy parameters maximize $\pi_x(d)$: the log-probability under $G_x$ of a graph with degree sequence $d$ is affine in its adjacency matrix with coefficients $\log(x_ix_j)$, hence equals its value at the maximum-entropy matrix $W$, which gives $\pi_x(d)=\pi(d)\exp\big(-\sum_{i<j}D(w^\star_{ij}\|w_{ij})\big)$ with $w^\star=W$ and $D$ the Bernoulli relative entropy of the following lemma.

\begin{lemma}[Pythagorean identity]\label{lem:pyth}
Let $d$ be strictly graphical and $D(y\|w)=y\log\frac yw+(1-y)\log\frac{1-y}{1-w}$. For every $Y\in\Psym(d)$,
\[
\sum_{i<j}D(y_{ij}\|w_{ij})\;=\;\gs(d)-\gsym(Y).
\]
Consequently $w_{ij}\ge\exp\!\big(-(\gs(d)-\gsym(Y)+\log2)/y_{ij}\big)$, and likewise for $1-w_{ij}$ with $1-y_{ij}$.
\end{lemma}
\begin{proof}
$\sum D(y\|w)=-\gsym(Y)-\sum_{i<j}\big[y_{ij}\log\tfrac{w_{ij}}{1-w_{ij}}+\log(1-w_{ij})\big]=-\gsym(Y)+\lambda\cdot d-\sum\log(1-w_{ij})$, using $\log\frac{w}{1-w}=-(\lambda_i+\lambda_j)$ and the margins of $Y$. The last two terms do not depend on $Y$; taking $Y=W$ shows they equal $\gs(d)$. For the bound, keep one term: $D(y\|w)\ge y\log\frac1w-\log2$, since $-(1-y)\log(1-w)\ge0$ and $y\log y+(1-y)\log(1-y)\ge-\log2$; all other terms are nonnegative.
\end{proof}

\begin{lemma}[Cut-norm concentration]\label{lem:cutconc}
Let $G$ have independent edges with probabilities $w_{ij}$ and $\sum_{i<j}w_{ij}\le M$. For $0<\varepsilon\le1$,
\[
\Pr\big[\|W^G-W^w\|_\square\ge\varepsilon\rho\big]\le2\cdot4^n\,e^{-\varepsilon^2M/3}.
\]
\end{lemma}
\begin{proof}
The cut norm of a step function on $\{I_i\}$ is attained at unions of intervals, by multilinearity. For $S,T\subseteq[n]$, $n^2\int_{S\times T}(W^G-W^w)=\sum_{i<j}c_{ij}(A_{ij}-w_{ij})$ with $c_{ij}\in\{0,1,2\}$; the summands are independent, bounded by $2$, with total variance at most $4M$. Bernstein's inequality at level $t=2\varepsilon M=\varepsilon\rho n^2$ gives probability at most $2\exp(-t^2/(2(4M+2t/3)))\le2e^{-\varepsilon^2M/3}$. Take a union bound over $4^n$ pairs.
\end{proof}

\begin{theorem}[Transfer of convergence]\label{thm:degconv}
Let $x=x^{(n)}$ satisfy $\sum_{i<j}w_{ij}\le M$, and suppose $M/n\to\infty$ and $-\log\pi_x(d)=o(M)$. Then $\|W^{G_d}/\rho-W^w/\rho\|_\square\to0$ in probability.
\end{theorem}
\begin{proof}
By \eqref{eq:transfer} and Lemma~\ref{lem:cutconc}, $\Pr[\|W^{G_d}-W^w\|_\square\ge\varepsilon\rho]\le2\exp(n\log4-\varepsilon^2M/3-\log\pi_x(d))\to0$.
\end{proof}

\begin{corollary}\label{cor:mainconv}
Let $d$ be strictly graphical. If $M/(n\log n\,\Lg)\to\infty$, then a uniformly random graph with degree sequence $d$ is, after normalization by $\rho$, cut-close to its maximum-entropy $\beta$-model with probability tending to $1$. Under Conjecture~\ref{conj:main} the same holds whenever $M/n\to\infty$.
\end{corollary}
\begin{proof}
By Corollary~\ref{cor:euler}, $-\log\pi(d)\le Cn\log n\,\Lg=o(M)$; under Conjecture~\ref{conj:main}, Theorem~\ref{thm:reduction}(2) gives $-\log\pi(d)\le(C+c_0/2)\,n\Lg=o(M)$. The maximum-entropy model has $\sum_{i<j}w_{ij}=M$. Apply Theorem~\ref{thm:degconv}.
\end{proof}

\subsection{An unconditional theorem in McKay's range}

\begin{theorem}[Sparse degree-constrained graphs converge]\label{thm:degsparse}
Suppose $M\to\infty$, $M/n\to\infty$ and $d_{\max}=o(M^{1/4})$. Then $\|W^{G_d}/\rho-\varphi_n\otimes\varphi_n\|_\square\to0$ in probability. In particular, if $\varphi_n\to\varphi$ in $L^1$, then $G_d$ converges in probability to the rank-one graphon $\varphi\otimes\varphi$ in the normalized cut metric.
\end{theorem}
\begin{proof}
Take $x_i=d_i/\sqrt{2M}$, so $w_{ij}=u_{ij}/(1+u_{ij})$ with $u_{ij}=d_id_j/(2M)\le d_{\max}^2/(2M)\to0$; then $\sum_{i<j}w_{ij}\le\sum_{i<j}u_{ij}\le M$. These are the parameters $q_{ij}$ of \cite[Thm.~2.25]{Meh14}.

\emph{Closeness to $\varphi_n\otimes\varphi_n$.} Off the diagonal blocks $\varphi_n\otimes\varphi_n=u_{ij}/\rho$ and $0\le(u_{ij}-w_{ij})/\rho\le u_{ij}^2/\rho$; on the diagonal blocks $W^w=0$. So $\|\varphi_n\otimes\varphi_n-W^w/\rho\|_1\le(\sum_{i,j}u_{ij}^2+\sum_iu_{ii})/(n^2\rho)\le(d_{\max}^2+d_{\max})/(2M)\to0$.

\emph{Lower bound on $\pi_x(d)$.} We have $\log\pi_x(d)=\log G(d)-\sum_{i<j}\log(1+u_{ij})+\sum_id_i\log x_i$. McKay's formula (Section~\ref{sec:sparse}) gives $G(d)=\frac{(2M)!}{M!2^M\prod d_i!}e^{-\lambda-\lambda^2+o(1)}$ with $\lambda\le d_{\max}/2$. Using $\log(1+u)\le u$, $\sum_{i<j}u_{ij}=M-\sum_id_i^2/(4M)$, $\log\frac{(2M)!}{M!2^M}\ge M\log(2M)-M-\log M-1$, $\log k!\le k\log k-k+\log k+1$ for $k\ge1$, and $\sum_{i:d_i\ge1}(1+\log d_i)\le n(1+\log\bar d)$, the terms $\pm\sum d_i\log d_i$ and $\pm M\log(2M)$ cancel, and since $\sum_id_i^2/(4M)-\lambda=\frac12$,
\[
-\log\pi_x(d)\le n(1+\log\bar d)+\log M+\frac{d_{\max}^2}4+1
\]
for large $n$. This is $o(M)$ because $\bar d\to\infty$ and $d_{\max}^2=o(M^{1/2})$.

Theorem~\ref{thm:degconv} and the triangle inequality give the first claim, and $\|\varphi_n\otimes\varphi_n-\varphi\otimes\varphi\|_1\le2\|\varphi_n-\varphi\|_1$ gives the second.
\end{proof}

\begin{remark}
\begin{enumerate}[leftmargin=2em,label=(\roman*)]
\item The same computation shows $\log G(d)=M\big(\log\frac{n^2}{2M}+1-2\int\varphi_n\log\varphi_n\big)+O(n\log\bar d+d_{\max}^2)$ in this range, which makes rigorous the expansion of \cite[Section~\ref{C-sec:program}]{Meh-open} there.
\item The limit $\varphi\otimes\varphi$ belongs to $L\log L$ iff $\varphi$ does, with $2\int\varphi\log\varphi-1\le\int\Psi(\varphi\otimes\varphi)\le2\int\varphi\log\varphi+2\log2$; combined with an Orlicz limit theory this places these graphs in the program of \cite[Section~\ref{C-sec:program}]{Meh-open}.
\item Beyond McKay's range, any bound $-\log\pi(d)=o(M)$ suffices for convergence to the $\beta$-model (Theorem~\ref{thm:degconv}); Corollary~\ref{cor:mainconv} gives it for $M\gg n\log n\,\Lg$. When $d_{\max}\gtrsim\sqrt M$ the top degrees saturate and the limit is the normalized maximum-entropy $\beta$-model rather than $\varphi\otimes\varphi$, although in the sparse regime the two have the same cut limit (\cite[Theorem~\ref{L-thm:betaphi}]{Meh-llogl}).
\item Local statistics: in \cite{Meh14} tree densities of uniformly random graphs were compared with the $\beta$-model, unconditionally for the relaxed constraint $|\deg_i-d_i|<d_i^a$ with $a>\frac12$, in McKay's range for exact degrees, and in general conditionally on \cite[Conj.~2.18]{Meh14}; by Theorem~\ref{thm:reduction}, that conjecture is equivalent (up to constants) to Conjecture~\ref{conj:main} with loss $e^{-O(n\log n)}$.
\end{enumerate}
\end{remark}

\section{Obstructions and consistency checks}

\subsection{Why Gurvits does not apply}
\begin{itemize}[leftmargin=2em]
\item \emph{Non-stability.} $\prod_{i<j}(1+x_ix_j)$ vanishes at $x_i=x_j=\sqrt{-1}$, so it is not real stable; its coefficient of $x^d$ is $G(d)$.
\item \emph{Parity.} Graphs need $\sum d_i$ even; digraphs do not. $d=(1,1,1)$ gives $D=2$, $G=0$.
\item \emph{No pointwise comparison.} Writing $D(d)=\sum_U\varepsilon(U)$ and $G(d)^2=\sum_U\kappa(U)$ over union multigraphs $U$ ($\varepsilon$ = Eulerian orientations, $\kappa$ = balanced red/blue edge colourings), an odd cycle has $\varepsilon=2$, $\kappa=0$. Any proof must compare the sums globally.
\end{itemize}

\subsection{Consistency checks}
\begin{center}
\renewcommand{\arraystretch}{1.25}
\begin{tabular}{@{}lllll@{}}
\toprule
Sequence & $G(d)$ & $D(d)$ & $G/\sqrt D$ & OK?\\
\midrule
perfect matching, $(1,\dots,1)$ & $(n-1)!!$ & $\sim n!/e$ & $\sim e^{1/2}(\pi n/2)^{-1/4}$ & yes (Thm.~\ref{thm:sparse} exactly)\\
$K_n$ & $1$ & $1$ & $1$ & yes\\
star $(n-1,1,\dots,1)$ & $1$ & $1$ & $1$ & yes\\
$d$-regular, $d$ fixed & Thm.~\ref{thm:sparse} & Thm.~\ref{thm:sparse} & $\sim e^{1/2}(\pi dn/2)^{-1/4}$ & yes\\
$(1,1,1)$ & $0$ & $2$ & $0$ & not graphical\\
\bottomrule
\end{tabular}
\end{center}
For threshold sequences $G(d)=1$, so the conjecture predicts $D(d)\le\exp(2Cn\Lg)$. If $\Psym(d)$ is a single point for threshold $d$, then $\gs=0$ and Corollary~\ref{cor:upper} gives $D(d)=1$; this has not been verified.

\end{document}